\documentclass[11pt]{article}
\usepackage{amscd}
\usepackage{amsfonts}
\usepackage{amsmath}
\usepackage{amssymb}
\usepackage{amsthm}
\usepackage{bbm}
\usepackage{CJK}
\usepackage{fancyhdr}
\usepackage{graphicx}
\usepackage{hyperref}
\usepackage{indentfirst}
\usepackage{latexsym}
\usepackage{mathrsfs}
\usepackage{xypic}

\usepackage[top=1in,bottom=1in,left=1.25in,right=1.25in]{geometry}
\def\<{\langle}
\def\>{\rangle}
\def\a{\alpha}
\def\b{\beta}

\def\c{\cdot}

\date{}
\begin{document}
\renewcommand{\baselinestretch}{1.2}
\renewcommand{\arraystretch}{1.0}
\title{\bf On the structure of split regular Hom-Lie Rinehart algebras}
\author{{\bf Shengxiang Wang$^{1}$,Xiaohui Zhang$^{2}$
        Shuangjian Guo$^{3}$\footnote
        { Corresponding author(Shuangjian Guo):~~shuangjianguo@126.com} }\\
1.~ School of Mathematics and Finance, Chuzhou University,\\
 Chuzhou 239000,  China \\
2.~  School of Mathematical Sciences, Qufu Normal University, \\Qufu 273165, China.\\
 3.~ School of Mathematics and Statistics, Guizhou University of\\ Finance and Economics,Guiyang 550025,China.}
 \maketitle
\begin{center}
\begin{minipage}{13.cm}

{\bf \begin{center} ABSTRACT \end{center}}
The aim of this paper is to study the structures of split regular Hom-Lie Rinehart algebras.
Let $(L,A)$  be a split regular Hom-Lie Rinehart algebra.
We first show that
 $L$ is of the form $L=U+\sum_{[\gamma]\in\Gamma/\thicksim}I_{[\gamma]}$ with $U$ a vector space complement in $H$
  and $I_{[\gamma]}$ are well described ideals of $L $ satisfying $I_{[\gamma]},I_{[\delta]}=0$ if $I_{[\gamma]}\neq I_{[\delta]}$.
Also, we discuss the weight spaces and decompositions of $A$ and present the relation between the decompositions of $L$ and $A$.
Finally, we consider the structures of tight split regular Hom-Lie Rinehart algebras.
 \smallskip

{\bf Key words}:  Hom-Lie Rinehart algebra; root space; weight  space; decomposition; simple ideal.
 \smallskip

 {\bf 2010 MSC:} 17A60; 17B22; 17B60; 17B65
 \end{minipage}
 \end{center}
 \normalsize\vskip0.5cm

\section{Introduction}
\def\theequation{\arabic{section}. \arabic{equation}}
\setcounter{equation} {0}

As an algebraic counterpart of Lie algebroids, the concept of Lie-Rinehart algebras was introduced
by Herz \cite{Herz}, Palais \cite{Palais} and Rinehart \cite{Rinehart} under the name of $(R,C)$-Lie algebra,
which play an important role in many branches of mathematics.
A Lie-Rinehart algebra $(L,A)$ can be viewed as a Lie  algebra $L$, which is simultaneously an $A$-module, where $A$ is an associative and
commutative algebra, in such a way that both structures are related in an appropriate way.
 For more details about the history and the developments  of Lie-Rinehart algebras, see
  \cite{Chemla, chen, Dokas, Huebschmann1999, Huebschmann2004, Rovi, Zhangtao}
and references cited therein.

There is a growing interest in twisted algebraic structures or Hom-algebraic structures defined for classical algebras
 and Lie algebroids as well.
 Hom-algebras were first introduced in the Lie algebra setting \cite{Hartwig} with motivation from physics
though its origin can be traced back in earlier literature such as \cite{Hu}.
In \cite{Makhlouf2008}, Makhlouf and Silvestrov introduced the definition of Hom-associative algebras,
where the associativity of a Hom-algebra is twisted by an endomorphism.
Recently,   Mandal  and Mishra introduced the notion of  Hom-Lie Rinehart algebras in \cite{Mandal},
which is a generalization of Lie-Rinehart algebras.

The class of the split algebras is specially related to addition quantum numbers, graded contractions and deformations.
For instance, for a physical system which displays a symmetry of $L$,
 it is interesting to know in detail the structure of the split decomposition because its roots can be seen as certain eigenvalues which are the additive quantum numbers characterizing the state of such system.
 Determining the structure of split algebras will become more and more meaningful in the area of research in mathematical physics.
Recently, in
 \cite{Albuquerque, Albuquerque2018, Aragon2016, Calderon2008, Calderon2009, Calderon2012, Calderon20121, Calderon2013, Cao2015, Cao2017, Cao2018},
 the structure of different classes of split algebras have been determined by the techniques of connections of roots.

In the present paper we introduce the class of split regular Hom-Lie Rinehart algebras as the natural extension of the one
of split Lie-Rinehart algebras and so of split regular Hom-Lie algebras, and study its tight structures based on some work in \cite{Albuquerque} and  \cite{Aragon2015}.
In section 2, we establish the preliminaries on split  regular Hom-Lie Rinehart algebras theory.
In sections 3 and 4, we  develop techniques of connections of roots and weights for split Hom-Lie Rinehart algebras respectively.
In section 5, we study the structures of tight split regular Hom-Lie Rinehart algebras.

\section{Preliminaries}
\def\theequation{\arabic{section}. \arabic{equation}}
\setcounter{equation} {0}

In this section, we start by recalling the definition of Hom-Lie Rinehart algebras,
then we introduce the notions of roots and weights of split Hom-Lie Rinehart algebras.
Throughout the paper, all algebraic systems are supposed to be over a field ${k}$.
\medskip

\noindent{\bf Definition 2.1.} (\cite{Okubo2002})
 Let  $A$ be an associative commutative algebra  and $\phi:A\rightarrow A$ an algebra endomorphism.
 A \emph{$\phi$-derivation} on $A$ is a  linear map $D:A\rightarrow A$  satisfying
\begin{eqnarray}
D(ab)=\phi(a)D(b)+D(a)\phi(b),
\end{eqnarray}
for all $a,b\in A$.
The set of all $\phi$-derivations on $A$ is denoted by $Der_{\phi}(A)$.
\medskip

\noindent{\bf Remark 2.2.}
 The triple $(Der_{\phi}(A),[\c,\c]_{\phi},\psi_{\phi})$ is a Hom-Lie algebra,
where the bracket $[\c,\c]_{\phi}$ and the twist map $\psi_{\phi}$ are defined by
 \begin{eqnarray*}
&&[D_1,D_2]_{\phi}=\phi\circ D_1\circ\phi^{-1}\circ D_2\circ\phi^{-1}-\phi\circ D_2\circ\phi^{-1}\circ D_1\circ\phi^{-1},\\
&&\psi_{\phi}(D)=\phi\circ D\circ\phi^{-1},
\end{eqnarray*}
for all $D_1,D_2\in Der_{\phi}(A)$.

Furthermore, $Der_{\phi}(A)$ has a natural $A$-module structure,
 where the $A$-module action is defined using the algebra multiplication in $A$.
These two structures are related by the following identities
 \begin{eqnarray*}
&&\psi_{\phi}(a\c D)=\phi(a)\c\psi_{\phi}(D),\\
&&[D_1,a\c D_2]_{\phi}=\phi(a)\c [D_1,D_2]_{\phi}+\psi_{\phi}(D_1)(a)\c\psi_{\phi}(D_2),
\end{eqnarray*}
for all $a\in A$ and $D_1,D_2\in Der_{\phi}(A)$.
\medskip

\noindent{\bf Definition 2.3.} (\cite{Mandal})
A \emph{Hom-Lie Rinehart algebra} over $(A,\phi)$ is a tuple $(A,L,[\c,\c],\phi,\psi,\rho)$,
where $A$ is an associative commutative algebra, $L$ is an $A$-module,
$[\c,\c]: L\times L\rightarrow L$ is a skew symmetric bilinear map,
the map $\phi: A\rightarrow A$ is an algebra homomorphism,
$\psi: L \rightarrow L$ is a linear map satisfying $\psi([x,y])=[\psi(x),\psi(y)]$,
 and $\rho: L \rightarrow Der_{\phi} A$ is a linear map,
satisfy the following conditions:

(1) the triple $(L,[\c,\c],\psi)$ is a Hom-Lie algebra.

(2) $\psi(a\c x)=\phi(a)\c\psi(x)$ for all $a\in A, x\in L.$

(3) $(\rho,\phi)$ is a representation of $(L,[\c,\c],\psi)$ on $A$.

(4) $\rho(a\c x)=\phi(a)\c\rho(x)$ for all $a\in A, x\in L.$

(5) $[x,a\c y]=\phi(a)[x,y]+\rho(x)(a)\psi(y)$ for all $a\in A, x,y\in L.$
\medskip

We denote it by $(L,A)$ or just by $L$ if there is not any possible confusion.
A Hom-Lie-Rinehart algebra $(A,L,[\c,\c],\phi,\psi,\rho)$ is said to be \emph{regular}
 if the map $\psi:L\rightarrow L$ is a bijective map and
  $\phi:A\rightarrow A $ is an algebra automorphism.
\medskip

\noindent{\bf Example 2.4.}
Any Hom-Lie algebra $(L,[\c,\c],\psi)$ over a $k$-module $L$ is a Hom-Lie Rinehart algebra over $A:= k$,
the algebra morphism $\phi=Id_k$  and the trivial action of $L$ on $k$.
\medskip

\noindent{\bf Example 2.5.}
Let $(L,A)$ be a Hom-Lie Rinehart algebra over an associative commutative $A$.
By Remark 2.2, $(Der_{\phi}(A),A,[\c,\c]_{\phi},\phi,\psi_{\phi},\rho_{\phi})$ is also a Hom-Lie Rinehart algebra with
$\rho_{\phi}=\psi_{\phi}$.
\medskip

\noindent{\bf Definition 2.6.}
 A \emph{subalgebra} $(S,A)$ of $(L,A)$, $S$ for short, is a Hom-Lie subalgebra of $L$ such that $AS\subset S$
 and satisfying that $S$ acts on $A$ via the composition $S\hookrightarrow L\rightarrow Der_{\phi} (A)$.
A subalgebra $(I,A)$, $I$ for short, of $L$ is called an \emph{ideal} if $I$ is a Hom-Lie ideal of $L$ such that $\rho(I)(A)L\subseteq I.$
\medskip

\noindent{\bf Definition 2.7.}
 A Hom-Lie Rinehart algebra $(L,A)$ is \emph{simple} if $[L,L]\neq 0,AA\neq 0,AL\neq 0$
  and its only ideals are ${0}, L $ and   the kernel of $\rho$.
\medskip

Throughout this paper $(L,A)$ will denote a regular  Hom-Lie Rinehart algebra over an arbitrary base field $k$ and that we
will denote by $N$ the set of all non-negative integers and by $Z $ the set of all integers.
\medskip

Next, we will introduce the class of split algebras in the framework of a  Hom-Lie Rinehart algebra  $(L,A)$.
We begin by recalling the definition of  a split regular Hom-Lie algebra.
\medskip

\noindent{\bf Definition 2.8.} (\cite{Aragon2015,Zhang2017})
Let $(L,[\c,\c],\psi)$ be a regular Hom-Lie algebra
and $H$ a maximal abelian subalgebra (MASA) of $L$.
We call $L$ is \emph{split} if  $L$ can be written as the direct sum
$L=H\oplus (\bigoplus_{\gamma\in\Gamma}L_{\alpha}),$
where
\begin{eqnarray}
L_{\gamma}=\{v_{\gamma}\in L|[h,v_{\gamma}]=\gamma(h)\psi(v_{\gamma}),~\forall h\in H\},
\end{eqnarray}
for a linear functional $\gamma\in H^{\ast}$ and $\Gamma=\{\gamma\in H^{\ast}\backslash\{0\}|L_\gamma\neq 0\}.$
$H$ is called  a \emph{splitting Cartan subalgebra} of $L$.
\medskip

Now let us return to a regular Hom-Lie Rinehart algebra  $(L,A)$.
Denote by $H$ a maximal abelian (in the sense $[H, H]=0$) subalgebra of $L$.
\medskip

\noindent{\bf Definition 2.9.}
 A split regular  Hom-Lie Rinehart algebra (with respect to a MASA $H$ of the regular Hom-Lie algebra $L$)
  is a \emph{regular  Hom-Lie Rinehart algebra} $(L,A)$ in which the Hom-Lie algebra $L$ contains a splitting Cartan
subalgebra $H$ and the algebra $A $ is a weight module (with respect to $H$) in the sense that
$A$ can be written as the direct sum
$A=A_0\oplus (\bigoplus_{\alpha\in \Lambda}A_{\alpha})$ with $\phi(A_\alpha)\subset A_\alpha$,
where
\begin{eqnarray}
A_{\alpha}=\{a_{\alpha}\in A|\rho(h)(a_{\alpha})=\alpha(h)\phi(a_{\alpha}),~\forall h\in H\},
\end{eqnarray}
for a linear functional $\alpha\in H^{\ast}$ and $\Lambda=\{\alpha\in H^{\ast}\backslash\{0\}|A_\gamma\neq 0\}$
 denotes the weights system of $A$. The linear subspace $ A_{\alpha}$ , for $\alpha\in \Lambda$, is called \emph{the weight space} of
$A$ associate to $\alpha$, the elements $\alpha\in \Lambda\cup \{0\}$ are called \emph{weights} of $A$.
\medskip

\noindent{\bf Remark 2.10.}
By \cite{Aragon2015, Zhang2017}, in a split regular  Hom-Lie Rinehart algebra $(L,A)$,
the  splitting Cartan subalgebra $H=L_0$,  that is,
$(L,A)$ has both a root  spaces decomposition $L=L_0\oplus (\bigoplus_{\gamma\in\Gamma}L_{\alpha})$
 and a weight  spaces decomposition $A=A_0\oplus (\bigoplus_{\alpha\in \Lambda}A_{\alpha}).$
 \medskip

\noindent{\bf Lemma 2.11.} For any $ \gamma,\xi\in\Gamma\cup \{0\}$  and $ \alpha,\beta\in\Lambda\cup \{0\}$,
 the following assertions hold.

(1)$L_0=H.$

(2) If $[L_{\gamma},L_{\xi}]\neq 0$, then
 $ \gamma\psi^{-1}+\xi\psi^{-1}\in\Gamma\cup \{0\}$ and $[L_{\gamma},L_{\xi}]\subset L_{\gamma\psi^{-1}+\xi\psi^{-1}}$.

(3) If $A_{\alpha}A_{\beta}\neq 0$, then
 $\alpha+\beta\in\Lambda\cup \{0\}$ and $A_{\alpha}A_{\beta}\subset A_{\alpha+\beta}$.

(4) If $A_{\alpha}L_{\gamma}\neq 0$, then
 $\alpha+\gamma\in\Gamma\cup \{0\}$ and $A_{\alpha}L_{\gamma}\subset L_{\alpha+\gamma}$.

(5) If $\rho(L_{\gamma})A_{\alpha}\neq 0$, then
$\alpha+\gamma\in\Lambda\cup \{0\}$ and $\rho(L_{\gamma})A_{\alpha}\subset A_{\alpha+\gamma}$.

\medskip

\noindent{\bf Proof.}
(1) Straightforward.

(2) For any $v_{\gamma}\in L_{\gamma},v_{\xi}\in L_{\xi}$ and $h\in H$, that is,
\begin{eqnarray*}
[h,v_{\gamma}]=\gamma(h)\psi(v_{\gamma}),[h,v_{\xi}]=\xi(h)\psi(v_{\xi}).
\end{eqnarray*}
By Hom-Jacobi identity, we have
\begin{eqnarray*}
[h,[v_{\gamma},v_{\xi}]]
&=&-[\psi(v_{\gamma}),[v_{\xi},\psi^{-1}(h)]]-[\psi(v_{\xi}),[\psi^{-1}(h),v_{\gamma}]]\\
&=&[\psi(v_{\gamma}),[\psi^{-1}(h),v_{\xi}]]+[[\psi^{-1}(h),v_{\gamma}],\psi(v_{\xi})]\\
&=&(\gamma\psi^{-1}+\xi\psi^{-1})(h)\psi[v_{\gamma},v_{\xi}].
\end{eqnarray*}
So $[v_{\gamma},v_{\xi}]\in L_{\gamma\psi^{-1}+\xi\psi^{-1}}$.

(3) For any $a_{\a}\in A_{\a}, a_{\b}\in A_{\b}$ and $h\in H$,
since $\rho(h)$ is a $\phi$-derivation on $A$, we have
\begin{eqnarray*}
\rho(h)(a_{\a}a_{\b})
&=&\phi(a_{\a})\rho(h)(a_{\b})+\rho(h)(a_{\a})\phi(a_{\b})\\
&=&\phi(a_{\a})\beta(h)\phi(a_{\b})+\alpha(h)\phi(a_{\a})\phi(a_{\b})\\
&=&(\alpha+\beta)(h)\phi(a_{\a}a_{\b}).
\end{eqnarray*}
Thus $a_{\a}a_{\b}\in A_{\alpha+\beta}$.

(4) For any $a_{\a}\in A_{\a},v_{\gamma}\in L_{\gamma}$  and $h\in H$,  we have
\begin{eqnarray*}
[h,a_{\a}v_{\gamma}]
&=&\phi(a_{\a})[h, v_{\gamma}]+\rho(h)(a_{\a})\psi(v_{\gamma})\\
&=&\gamma(h)\phi(a_{\a})\psi(v_{\gamma})+\alpha(h)\phi(a_{\a})\psi(v_{\gamma})\\
&=&(\gamma+\alpha)(h)\psi(a_{\a}v_{\gamma}).
\end{eqnarray*}
Therefore $a_{\a}v_{\gamma}\in L_{\alpha+\gamma}$.

(5) For any $a_{\a}\in A_{\a},v_{\gamma}\in A_{\gamma}$  and $h\in H$,
since  $(\rho,\phi)$ is a representation of $(L,[\c,\c],\psi)$ on $A$, we have
\begin{eqnarray*}
\rho([h,v_{\gamma}])\phi(a_{\a})
=\rho(\psi(h))\rho(v_{\gamma})(a_{\a})-\rho(\psi(v_{\gamma}))\rho(h)(a_{\a}).
\end{eqnarray*}
Because  $a_{\a}\in A_{\a},v_{\gamma}\in L_{\gamma}$, we have
\begin{eqnarray*}
&&\rho([h,v_{\gamma}])\phi(a_{\a})
=\rho(\gamma(h)v_{\gamma})\phi(a_{\a})
=\gamma(h)\phi \rho(v_{\gamma})\phi(a_{\a}),\\
&&\rho(\psi(v_{\gamma}))\rho(h)(a_{\a})
=\rho(\psi(v_{\gamma}))\alpha(h)\phi(a_{\a})
=\alpha(h)\phi \rho(v_{\gamma})\phi(a_{\a}).
\end{eqnarray*}
It follows that
\begin{eqnarray*}
(\alpha(h)+\gamma(h))\phi \rho(v_{\gamma})\phi(a_{\a})
=\rho(\psi(h))\rho(v_{\gamma})(a_{\a})
=\phi\rho(h)\phi^{-1}\rho(v_{\gamma})(a_{\a}).
\end{eqnarray*}
The second equality holds since $\rho$ is a homomorphism of Hom-Lie algebras.
Therefore, $\rho(h)(\rho(v_{\gamma})(a_{\a}))=(\alpha(h)+\gamma(h))\phi(\rho(v_{\gamma})\phi(a_{\a})).$
That is,  $\rho(v_{\gamma})(a_{\a})\in A_{\alpha+\gamma}$.
$\hfill \Box$

\section{Connections of roots}
\def\theequation{\arabic{section}. \arabic{equation}}
\setcounter{equation} {0}

In the following, $(L,A)$ denotes a split regular Hom-Lie Rinehart algebra
which has  a root  spaces decomposition  and a weight  spaces decomposition:
\begin{eqnarray}
L=L_0\oplus (\bigoplus_{\gamma\in\Gamma}L_{\alpha}),~
A=A_0\oplus (\bigoplus_{\alpha\in \Lambda}A_{\alpha}).
\end{eqnarray}
Given a linear functional $\gamma: H\rightarrow  k$, we denote by $-\gamma: H\rightarrow  k$ the element
in $H^{\ast}$ defined by $(-\gamma)(h):=-\gamma(h)$ for all $h\in  H.$ We also denote $-\Gamma:=\{-\gamma|\gamma\in \Gamma\}$.
In a similar way we can define $-\Lambda:=\{-\alpha|\alpha\in\Lambda\}$.
Finally, we denote $\pm\Gamma:=\Gamma\cup -\Gamma $ and $\pm\Lambda:=\Lambda\cup -\Lambda.$
\medskip

\noindent{\bf Definition 3.1.}
Let $\gamma,\xi\in \Gamma$,
 we say that $\gamma$ and $\xi$ are \emph{connected} if

 $\bullet$ Either $\xi=\epsilon\gamma\psi^{z}$ for some $z\in Z$ and $\epsilon\in \{1,-1\}$.

  $\bullet$ Either there exists a family $\{\zeta_1,\zeta_2,\cdots,\zeta_n\}\subset\pm\Lambda\cup\pm\Gamma$,
     with $n\geq 2$, such that

(1) $\zeta_1\in\{\gamma\psi^{k}|k\in Z\}$.

(2) $\zeta_1\psi^{-1}+\zeta_2\psi^{-1}\in\pm\Gamma$,

    $~~~~~\zeta_1\psi^{-2}+\zeta_2\psi^{-2}+\zeta_3\psi^{-1}\in\pm\Gamma$,

    $~~~~~\zeta_1\psi^{-3}+\zeta_2\psi^{-3}+\zeta_3\psi^{-2}+\zeta_4\psi^{-1}\in\pm\Gamma$,

    $~~~~~~~~~~\cdots\cdots\cdots$

    $~~~~~\zeta_1\psi^{-i}+\zeta_2\psi^{-i}+\zeta_3\psi^{-i+1}+\cdots+\zeta_{i+1}\psi^{-1}\in\pm\Gamma$,

    $~~~~~~~~~~\cdots\cdots\cdots$

    $~~~~~\zeta_1\psi^{-n+2}+\zeta_2\psi^{-n+2}+\zeta_3\psi^{-n+3}+\cdots+\zeta_{n-1}\psi^{-1}\in\pm\Gamma$.

(3) $\zeta_1\psi^{-n+1}+\zeta_2\psi^{-n+1}+\zeta_3\psi^{-n+2}+\cdots+\zeta_{n}\psi^{-1}\in\{\pm\xi\psi^{-m}|m\in Z\}$.

We will also say that $\{\zeta_1,\zeta_2,\cdots,\zeta_n\}$ is a \emph{connection from $\gamma$ to $\xi$}.
\medskip

\noindent{\bf Remark 3.2.}
Assume that $\gamma$ and $\xi$ are connected.
Then $\gamma\psi^{z_1}$  is  connected to $\xi\psi^{z_2}$ for any $z_1,z_2\in Z$,
and also to $-\xi\psi^{z_2}$ in case $-\xi\psi^{z_2}\in \Gamma$.

\medskip

\noindent{\bf Proposition 3.3.}
The relation $\sim$ in $\Gamma$ is an equivalence relation, where $\gamma\sim\xi$ if and only if $\gamma$  is  connected to $\xi$.
 \medskip

\noindent{\bf Proof.}
Similar to Proposition 2.4 in \cite{Aragon2015}.
$\hfill \Box$
 \medskip

\noindent{\bf Remark 3.4.}
Let $\xi,\gamma\in\Gamma$ such that $\gamma\sim\xi$,
$\eta\in\Lambda\cap\Gamma$ such that $\gamma+\eta\in\Gamma$.
Since $\{\gamma,\eta\}$ is a collection from $\gamma$ to $\gamma+\eta$,
we get $\xi\sim\gamma+\eta.$
 \medskip

By Proposition 3.3, we can consider the quotient set
$$\Gamma/\sim:=\{[\gamma]|\gamma\in \Gamma\},$$
where $[\gamma]$ denotes the set of nonzero roots of $L$ which are connected to $\gamma$.
In the following we will  associate an  adequate  ideal $I_{[\gamma]}$ to any $[\gamma]$.
 For a fixed $\gamma\in \Gamma$, we define
\begin{eqnarray}
L_{0,[\gamma]}:=(\sum_{\xi\in [\gamma],-\xi\in\Lambda}A_{-\xi}L_{\xi})+(\sum_{\xi\in [\gamma]}[L_{-\xi},L_{\xi}])\subset L_0,~~
 L_{[\gamma]}:=\bigoplus_{\xi\in [\gamma]}L_{\xi}.
 \end{eqnarray}
Then we denote by $I_{[\gamma]}$ the direct sum of the two subspaces above, that is,
 \begin{eqnarray}
 I_{[\gamma]}:=L_{0,[\gamma]}\oplus L_{[\gamma]}.
 \end{eqnarray}

\noindent{\bf Proposition 3.5.}
For any $\gamma\in \Gamma$, the  following assertions hold.

(1) $[I_{[\gamma]},I_{[\gamma]}]\subset I_{[\gamma]}$.

(2) $\psi(I_{[\gamma]})=I_{[\gamma]}$.

(3) $AI_{[\gamma]}\subset I_{[\gamma]}$.

(4) $\rho(I_{[\gamma]})(A)L\subset I_{[\gamma]}$.
\medskip

\noindent{\bf Proof}
(1) Since $L_{0,[\gamma]}\subset L_0=H$, we have $[L_{0,[\gamma]},L_{0,[\gamma]}]=0$  and therefore
\begin{eqnarray}
[I_{[\gamma]},I_{[\gamma]}]
=[L_{0,[\gamma]}\oplus L_\gamma,L_{0,[\gamma]}\oplus L_\gamma]
\subset[L_{0,[\gamma]}, L_\gamma]+[ L_\gamma,L_{0,[\gamma]}]+[L_\gamma,L_\gamma].
\end{eqnarray}

For any $\delta\in [\gamma]$, by Eq. (3.2) and Lemma 2.11-(2),
we have $[L_{0,[\gamma]},L_\delta]\subset L_{\delta\psi^{-1}}$, where $\delta\psi^{-1}\in [\gamma]$.
It follows that $[L_{0,[\gamma]},L_{[\gamma]}]\subset L_{[\gamma]}$.
Similarly, we have  $[L_{[\gamma]}, L_{0,[\gamma]}]\subset L_{[\gamma]}$.

We now consider the expression $[L_\gamma,L_\gamma]$ in Eq. (3.4).
For this, we take any $\delta,\eta\in [\gamma]$ satisfying $[L_\delta,L_\eta]\neq 0$,
by Lemma 2.11-(2), $[L_\delta,L_\eta]\subset L_{\delta\psi^{-1}+\eta\psi^{-1}}$.

$\bullet$ If $\delta+\eta=0$, then  $[L_\delta,L_\eta]\subset L_{0,[\gamma]}$.

$\bullet$ If $\delta+\eta\neq 0$, by Lemma 2.11-(2), we have $\delta\psi^{-1}+\eta\psi^{-1}\in \Gamma$ and therefore
$\{\delta,\eta\}$ is a collection from $\delta$ to $\delta\psi^{-1}+\eta\psi^{-1}$.
By Proposition 3.3,  $\delta\psi^{-1}+\eta\psi^{-1}\in [\gamma]$.
Thus $[L_\delta,L_\eta]\subset L_{\delta\psi^{-1}+\eta\psi^{-1}}\subset L_{[\gamma]}$.

So we get $[I_{[\gamma]},I_{[\gamma]}]\subset I_{[\gamma]}$.
\medskip

(2) For any $\gamma\in \Gamma$, we first claim that
\begin{eqnarray}
\psi(L_\gamma)\subset L_{\gamma\psi^{-1}},~ \psi^{-1}(L_\gamma)\subset L_{\gamma\psi}.
\end{eqnarray}
In fact, for any $h\in H$ and $v_{\gamma}\in  L_\gamma$, we have
\begin{eqnarray*}
[h,\psi(v_{\gamma})]
=\psi([\psi^{-1}(h),v_{\gamma}])
=\gamma(\psi^{-1}(h))\psi(\psi(v_{\gamma})).
\end{eqnarray*}
It follows that $\psi(v_{\gamma})\in L_{\gamma\psi^{-1}}$. Thus $\psi(L_\gamma)\subset L_{\gamma\psi^{-1}}$.
Similarly, one may check that $\psi^{-1}(L_\gamma)\subset L_{\gamma\psi}$.
By Remark 3.2, $\gamma\psi^{-1}\in [\gamma], \gamma\psi\in [\gamma]$,
then we have $\psi(I_{[\gamma]})=I_{[\gamma]}$, as desired.
\medskip

(3) By Eqs. (3.1) and (3.3), we have
\begin{eqnarray*}
AI_{[\gamma]}=(A_0\oplus (\bigoplus_{\alpha\in \Lambda}A_{\alpha}))
((\sum_{\xi\in [\gamma],-\xi\in\Lambda}A_{-\xi}L_{\xi})+(\sum_{\xi\in [\gamma]}[L_{-\xi},L_{\xi}])+(\bigoplus_{\xi\in [\gamma]}L_{\xi})).
\end{eqnarray*}
We have to consider six cases:

Case 1: $A_0(A_{-\xi}L_{\xi})\subset L_{0,[\gamma]}.$
In fact, for any $\xi\in [\gamma]$ and $-\xi\in\Lambda$, by Lemma 2.11-(2-3) and the $A$-module structure on $L$, we have
\begin{eqnarray*}
A_0(A_{-\xi}L_{\xi})
=(A_0A_{-\xi})L_{\xi}
\subset A_{-\xi}L_{\xi}
\subset L_{0,[\gamma]}.
\end{eqnarray*}

Case 2: $A_0[L_{-\xi},L_{\xi}]\subset L_{0,[\gamma]}.$
In fact, for any $\xi\in [\gamma]$ and $ a_0\in A_0, v_{-\xi}\in L_{-\xi}, v_{\xi}\in L_{\xi}$, by  Definition 2.3, we have
\begin{eqnarray*}
a_0[v_{-\xi},v_{\xi}]
=[v_{-\xi},\phi^{-1}(a_0)v_{\xi}]-\rho(v_{-\xi})(\phi^{-1}(a_0))\psi(v_{\xi}).
\end{eqnarray*}
Since $\phi^{-1}(a_0)\in A_0$, we have $\phi^{-1}(a_0)v_{\xi}\in L_\xi$
and therefore $[v_{-\xi},\phi^{-1}(a_0)v_{\xi}]\in [L_{-\xi},L_{\xi}].$
Also, we can get $\rho(v_{-\xi})(\phi^{-1}(a_0))\in A_{-\xi}.$
If  $A_{-\xi}\neq 0$ (otherwise is trivial), $-\xi\in \Lambda $.
By Lemma 2.11-(5), we have
$
\rho(v_{-\xi})(\phi^{-1}(a_0))\psi(v_{\xi})
\subset  A_{-\xi}L_{\xi\psi^{-1}},
$
where $\xi\psi^{-1}\in [\gamma]$. So $A_0[L_{-\xi},L_{\xi}]\subset L_{0,[\gamma]}.$

Case 3: $A_0L_{\xi}\subset L_{\xi}\subset L_{[\gamma]}$. Straightforward from Lemma 2.11-(4).

Case 4: $A_\alpha(A_{-\xi}L_{\xi})\subset L_{[\gamma]}.$
In fact, for any $\alpha\in\Lambda$, $\xi\in [\gamma]$ and $-\xi\in\Lambda$, by Lemma 2.11-(2-3) and  the $A$-module structure on $L$, we have
\begin{eqnarray*}
A_\alpha(A_{-\xi}L_{\xi})
=(A_\alpha A_{-\xi})L_{\xi}
\subset A_{\alpha-\xi}L_{\xi}
\subset L_{\alpha}.
\end{eqnarray*}
If  $L_{\alpha}\neq 0$ (otherwise is trivial), then $\alpha\in \Gamma $. By Remark 3.4, $\alpha\in [\gamma]. $
Thus we have $A_\alpha(A_{-\xi}L_{\xi})\subset L_{[\gamma]}.$

Case 5: $A_\alpha[L_{-\xi},L_{\xi}]\subset L_{[\gamma]}.$
In fact, for any $\alpha\in\Lambda, \xi\in [\gamma]$ and $ a_\alpha\in A_\alpha, v_{-\xi}\in L_{-\xi}, v_{\xi}\in L_{\xi}$,
 by  Definition 2.3, we have
\begin{eqnarray*}
a_\alpha[v_{-\xi},v_{\xi}]
=[v_{-\xi},\phi^{-1}(a_\alpha)v_{\xi}]-\rho(v_{-\xi})(\phi^{-1}(a_\alpha))\psi(v_{\xi}).
\end{eqnarray*}

Since $\phi^{-1}(a_\alpha)\in A_\alpha$, we have $\phi^{-1}(a_\alpha)v_{\xi}\in L_{\alpha+\xi}$.
If $\rho(v_{-\xi})(\phi^{-1}(a_\alpha))\neq 0$ (otherwise is trivial), $\alpha+\xi\in \Gamma$.
If $[v_{-\xi},\phi^{-1}(a_\alpha)v_{\xi}]\neq 0$ (otherwise is trivial), then $\alpha\psi^{-1}\in \Gamma$ and
$[v_{-\xi},\phi^{-1}(a_\alpha)v_{\xi}]\in [L_{-\xi},L_{\alpha+\xi}]\subset L_{\alpha\psi^{-1}}.$
By Remark 3.3,  $\alpha+\xi\in [\gamma]$,
then $\alpha\sim \alpha+\xi$, thus $\alpha\in [\gamma]$ and therefore $\alpha\psi^{-1}\in [\gamma]$.
So $[v_{-\xi},\phi^{-1}(a_\alpha)v_{\xi}]\in  L_{[\gamma]}.$

If  $\rho(v_{-\xi})(\phi^{-1}(a_\alpha))\neq 0$ (otherwise is trivial),
then  $\alpha-\xi\in \Lambda$ and $\rho(v_{-\xi})(\phi^{-1}(a_\alpha))\in A_{\alpha-\xi}.$
Similar to the discussion above,  $\alpha\in [\gamma]$.
If  $\rho(v_{-\xi})(\phi^{-1}(a_\alpha))\psi(v_{\xi})\neq 0$ (otherwise is trivial),
by Lemma 2.11-(4), we have $\alpha-\xi+\xi\psi^{-1}\in \Gamma$ and
$
\rho(v_{-\xi})(\phi^{-1}(a_\alpha))\psi(v_{\xi})
\in  L_{\alpha-\xi+\xi\psi^{-1}},
$
where $\xi,\xi\psi^{-1}\in [\gamma]$. So $A_0[L_{-\xi},L_{\xi}]\subset L_{[\gamma]}.$

Case 6:  $A_\alpha L_{\xi}\subset L_{[\gamma]}.$
In fact, for any $\alpha\in\Lambda, \xi\in [\gamma]$, we have $A_\alpha L_{\xi}\subset L_{\xi+\alpha}.$
Using again Remark 3.3, $\xi+\alpha\in [\gamma]$. So $A_\alpha L_{\xi}\subset L_{[\gamma]}.$

According to the six cases above, we prove $AI_{[\gamma]}\subset I_{[\gamma]}$, as required.

(4) By the definition 2.3, we have
$
\rho(I_{[\gamma]})(A)L
\subset [I_{[\gamma]},AL]+\phi(A)[I_{[\gamma]},L].
$
By the conclusion (3) above, we get $\rho(I_{[\gamma]})(A)L\subset I_{[\gamma]}.$ $\hfill \Box$
\medskip

\noindent{\bf Proposition 3.6.}
If $[\gamma],[\delta]\in \Gamma/\sim$ such that $[\gamma]\neq[\delta]$,  then $[I_{[\gamma]},I_{[\delta]}]=0$.
\medskip

\noindent{\bf Proof.}
Since $I_{[\gamma]}=L_{0,[\gamma]}\oplus L_\gamma$ and $L_{0,[\gamma]}\subset L_0$, we have
\begin{eqnarray}
[I_{[\gamma]},I_{[\gamma]}]
 =[L_{0,[\gamma]}\oplus L_\gamma,L_{0,[\delta]}\oplus L_\delta]
\subset [L_{0,[\gamma]},L_\delta]+[L_\gamma,L_{0,[\delta]}]+[L_\gamma,L_\delta].
\end{eqnarray}

We first claim that $[L_\gamma,L_\delta]=0.$
Suppose there exist $\gamma_1\in[\gamma],\delta_1\in[\delta]$ such that $[L_{\gamma_1},L_{\delta_1}]\neq 0,$
by Lemma 2.11-(2), we get
$[L_{\gamma_1},L_{\delta_1}]\subset L_{\gamma_1\psi^{-1}+\delta_1\psi^{-1}}=L_{(\gamma_1+\delta_1)\psi^{-1}}$
and $\gamma_1+\delta_1\neq 0.$
Thus $\gamma_1+\delta_1\in \Gamma$.
Since $\gamma_1\in[\gamma],\delta_1\in[\delta]$, by remark 3.3,
we have $\gamma\sim \gamma_1+\delta_1, \delta\sim \gamma_1+\delta_1$
and therefore $\gamma\sim\delta$, contradiction. So $[L_\gamma,L_\delta]=0.$

Second, we claim that $[L_{0,[\gamma]},L_\delta]=[L_\gamma,L_{0,[\delta]}]=0$.
By the definition of $L_{0,[\gamma]}, L_\delta$, we have
$$[L_{0,[\gamma]},L_\delta]=[(\sum_{\gamma_1\in [\gamma],-\gamma_1\in\Lambda}A_{-\gamma_1}L_{\gamma_1})
+(\sum_{\gamma_{1}\in [\gamma]}[L_{-\gamma_1},L_{\gamma_1}]),L_\delta].$$

For the expression $[L_{0,[\gamma]},L_\delta]$ in Eq. (3.6), we take any  $\delta_1\in [\delta]$, by the Hom-Jacobi identity,
$$[L_{-\gamma_1},L_{\gamma_1}],L_{\delta_1}]=[[L_{\gamma_1},L_{\delta_{1}\psi}],L_{-\gamma_{1}\psi^{-1}}]
+[[L_{\delta_{1}\psi},L_{-\gamma_1}],L_{\gamma_{1}\psi^{-1}}].$$
Since $[L_{\gamma_1},L_{\delta_{1}\psi}],L_{-\gamma_{1}\psi^{-1}}]=[L_{\delta_{1}\psi},L_{-\gamma_1}]=0$,
we have $[L_{0,[\gamma]},L_\delta]=0$.

For the expression  $[A_{-\gamma_1}L_{\gamma_1},L_\delta]$ in Eq. (3.6), suppose there exists $\delta_1\in [\delta]$ such that $[A_{-\gamma_1}L_{\gamma_1},L_\delta]\neq 0$.
 By Definition 2.3, we have
$$[A_{-\gamma_1}L_{\gamma_1},L_{\delta_1}]
=[L_{\delta_1},A_{-\gamma_1}L_{\gamma_1}]
\subset \phi(A_{-\gamma_1})[L_{\delta_1},L_{\gamma_1}]+\rho(L_{\delta_1})(A_{-\gamma_1})L_{\delta_{1}\psi^{-1}}.$$
By the discussion above, we get $[L_{\delta_1},L_{\gamma_1}]=0$.
Since $[A_{-\gamma_1}L_{\gamma_1},L_\delta]\neq 0$,
it follows that $0 \neq\rho(L_{\delta_1})(A_{-\gamma_1})L_{\delta_{1}\psi^{-1}}\subset A_{\delta_{1}-\gamma_{1}}L_{\delta_{1}\psi^{-1}}$.
Thus $A_{\delta_{1}-\gamma_{1}}\neq 0$ and $\delta_{1}-\gamma_{1}\in \Lambda\cup\{0\}$.
By Remark 3.3, $\delta_{1}\sim\gamma_{1}$, a contradiction.
So  $[A_{-\gamma_1}L_{\gamma_1},L_\delta]=0$.

Therefore, we have $[L_{0,[\gamma]},L_\delta]=0$.
Similarly, one may check that $[L_\gamma,L_{0,[\delta]}]=0$.
The proof is finished.
$\hfill \Box$
\medskip

\noindent{\bf Theorem 3.7.}
The following assertions hold.

(1) For any $[\gamma]\in \Gamma/\sim$, the linear space $I_{[\gamma]}=L_{0,[\gamma]}\oplus L_\gamma$
associated to $[\gamma]$ is an ideal of $L$.

(2) If $L$ is simple, then there is a connection from $\gamma$ to $\xi$ for any $\gamma,\xi\in \Gamma$.
Moreover,
\begin{eqnarray}
H=(\sum_{\gamma\in \Gamma,-\gamma\in\Lambda}A_{-\gamma}L_{\gamma})+(\sum_{\gamma\in \Gamma}[L_{-\gamma},L_{\gamma}]).
\end{eqnarray}

\noindent{\bf Proof.}
(1) Since $[I_{[\gamma]},H]=[I_{[\gamma]},L_0]\subset L_{[\gamma]}\subset I_{[\gamma]}$,
taking into account Propositions 3.5 and 3.6,  we have
$$[I_{[\gamma]},L]=[I_{[\gamma]}, H\oplus (\bigoplus_{\xi\in[\gamma]}L_{\xi})\oplus (\bigoplus_{\delta\notin[\gamma]}L_{\delta})\subset I_{[\gamma]}.$$
Since $\psi(L_\gamma)\subset L_{\gamma\psi^{-1}}$  and $\psi^{-1}(L_\gamma)\subset L_{\gamma\psi}$,
by Remark 3.2, we have $\psi(I_{[\gamma]})=I_{[\gamma]}$.
By Proposition 3.5, $I_{[\gamma]}$ is an ideal of $L$.

(2)  The simplicity of $L$ implies $I_{[\gamma]}=Ker\rho$ or $I_{[\gamma]}=L$ for any $\gamma\in \Gamma$.
 If there exists $\gamma\in \Gamma$ such that $I_{[\gamma]}=L$, then $[\gamma]=\Gamma$.
 Otherwise, if $I_{[\gamma]}=Ker\rho$ for all $\gamma\in \Gamma$  then $[\gamma]=[\xi]$ for any
 $\gamma,\xi\in \Gamma$ and again $[\gamma]=\Gamma$.
 Therefore in any case $L$ has all its nonzero roots connected and
$H=(\sum_{\gamma\in \Gamma,-\gamma\in\Lambda}A_{-\gamma}L_{\gamma})+(\sum_{\gamma\in \Gamma}[L_{-\gamma},L_{\gamma}]).$
$\hfill \Box$
\medskip

\noindent{\bf Theorem 3.8.}
For a vector space complement $U$ of
$(\sum_{\gamma\in \Gamma,-\gamma\in\Lambda}A_{-\gamma}L_{\gamma})+(\sum_{\gamma\in \Gamma}[L_{-\gamma},L_{\gamma}])$ in $H$,
we have
$$L=U+\sum_{[\gamma]\in\Gamma/\thicksim}I_{[\gamma]},$$
where $I_{[\gamma]}$  is one of the ideals of $L$ described in Theorem 3.7-(1),
satisfying $[I_{[\gamma]},I_{[\delta]}]=0$, whenever $[\gamma]\neq [\delta].$
\medskip

\noindent{\bf Proof.}
 Each $I_{[\gamma]}$ is well defined and, by Theorem 3.7-(1), an   ideal of L.
It is clear that
$$L=U+\sum_{\gamma\in\Gamma}L_{[\gamma]}=U+\sum_{[\gamma]\in\Gamma/\thicksim}I_{[\gamma]},$$
By Proposition 3.6, we get  $[I_{[\gamma]},I_{[\delta]}]=0$ if $[\gamma]\neq [\delta].$
$\hfill \Box$
\medskip

Let us denote by $Z(L)$ the center of $L$, that is,
$Z(L):=\{v\in L| [v, L]=0, \rho(v)=0\}.$
\medskip

\noindent{\bf Corollary 3.9.}
 If $Z(L)=0$ and $H=(\sum_{\gamma\in \Gamma,-\gamma\in\Lambda}A_{-\gamma}L_{\gamma})+(\sum_{\gamma\in \Gamma}[L_{-\gamma},L_{\gamma}])$,
  then $L$ is the direct sum of the ideals given in Theorem 3.6,
$$L=\bigoplus_{[\gamma]\in\Gamma/\thicksim}I_{[\gamma]}.$$

\noindent{\bf Proof.}
Straightforward from Theorem 3.8.
$\hfill \Box$

\section{Decompositions of  $A$}
\def\theequation{\arabic{section}. \arabic{equation}}
\setcounter{equation} {0}

Corresponding to the   decompositions of $L$ in the split regular  Hom-Lie Rinehart algebra $(L,A)$,
we will discuss the weight  spaces and decompositions of $A$ in the following.
\medskip

\noindent{\bf Definition 4.1.}
Let $\alpha,\beta\in \Lambda$
 we say that $\alpha$ and $\beta$ are \emph{connected} if

 $\bullet$ either $\beta=\varepsilon\alpha$ for some   $\varepsilon\in \{1,-1\}$;

  $\bullet$ either there exists a family $\{\sigma_1,\sigma_2,\cdots,\sigma_n\}\subset\pm\Lambda\cup\pm\Gamma$,
     with $n\geq 2$, such that

(1) $\sigma_1=\alpha$.

(2) $\sigma_1+\sigma_2\in\pm\Lambda\cup\pm\Gamma$,

    $~~~~~\sigma_1+\sigma_2+\sigma_3\in\pm\Lambda\cup\pm\Gamma$,

    $~~~~~~~~~~\cdots\cdots\cdots$

    $~~~~~\sigma_1+\sigma_2\cdots+\sigma_{n-1}\in\pm\Lambda\cup\pm\Gamma$,

(3) $\sigma_1+\sigma_2\cdots+\sigma_{n}\in\{\beta,-\beta\}$.

We will also say that $\{\sigma_1,\sigma_2,\cdots,\sigma_n\}$ is a \emph{connection} from $\alpha$ to $\beta$.
\medskip

\noindent{\bf Proposition 4.2.}
The relation $\approx$ in $\Lambda$ is an equivalence relation, where $\alpha\approx\beta$ if and only if
 $\alpha$  is  connected to $\beta$.
 \medskip

\noindent{\bf Proof.} Straightforward.
$\hfill \Box$
 \medskip

\noindent{\bf Remark 4.3.}
Let $\alpha,\beta\in\Lambda$ such that $\alpha\approx\beta$,
$\mu\in\Lambda\cap\Gamma$ such that $\beta+\mu\in\Lambda$.
Since $\{\beta,\mu\}$ is a collection from $\beta$ to $\beta+\mu$,
we get $\alpha\approx\beta+\mu.$
 \medskip

By Proposition 4.2, we can consider the quotient set
$$\Lambda/\approx:=\{[\alpha]|\alpha\in \Lambda\},$$
where $[\alpha]$ denotes the set of nonzero weights of $A$ which are connected to $\alpha$.
In the following we will  associate an  adequate  ideal $\mathcal{A}_{[\gamma]}$ to any $[\alpha]$.
 For a fixed $\alpha\in \Lambda$, we define
\begin{eqnarray}
A_{0,[\gamma]}:=(\sum_{\beta\in [\alpha],-\beta\in\Lambda}\rho(L_{-\beta})(A_{\beta})
                    +(\sum_{\beta\in [\alpha]}A_{-\beta},A_{\beta})\subset A_0,~~
 A_{[\alpha]}:=\bigoplus_{\beta\in [\alpha]}A_{\beta}.
 \end{eqnarray}
Then we denote by $\mathcal{A}_{[\alpha]}$ the direct sum of the two subspaces above,
 \begin{eqnarray}
 \mathcal{A}_{[\alpha]}:=A_{0,[\alpha]}\oplus A_{[\alpha]}.
 \end{eqnarray}

\noindent{\bf Proposition 4.4.}
For any $\alpha,\beta\in \Lambda$, the following assertions hold.

(1) $\mathcal{A}_{[\alpha]}\mathcal{A}_{[\alpha]}\subset \mathcal{A}_{[\alpha]}$.

 (2) If $[\alpha]\neq[\beta]$, then $\mathcal{A}_{[\alpha]}\mathcal{A}_{[\beta]}=0$.
 \medskip

\noindent{\bf Proof.}
We only prove (1) and similar for (2).
By Eq. (4.2)  and the commutativity of $A$, we have
\begin{eqnarray}
\mathcal{A}_{[\gamma]}\mathcal{A}_{[\gamma]}
=(A_{0,[\alpha]}\oplus A_[\alpha])(A_{0,[\alpha]}\oplus A_{[\alpha]})
\subset A_{0,[\alpha]}A_{0,[\alpha]}+A_{0,[\alpha]}A_{[\alpha]}+A_{[\alpha]}A_{[\alpha]}.
\end{eqnarray}
We first consider the expression $A_{0,[\alpha]}A_{0,[\alpha]}$ in Eq. (4.3).
Let  $\beta,\gamma\in[\alpha]$ satisfying
\begin{eqnarray*}
&&(\rho(L_{-\beta})(A_{\beta})+A_{-\beta}A_{\beta})(\rho(L_{-\gamma})(A_{\gamma})+A_{-\gamma}A_{\gamma})\\
&=&(\rho(L_{-\beta})(A_{\beta}))(\rho(L_{-\gamma})(A_{\gamma}))
   +(\rho(L_{-\beta})(A_{\beta}))(A_{-\gamma}A_{\gamma})\\
   &&+(A_{-\beta}A_{\beta})(\rho(L_{-\gamma})(A_{\gamma}))
   +(A_{-\beta}A_{\beta})(A_{-\gamma},A_{\gamma})\neq 0.
\end{eqnarray*}
There are four cases we have to discuss:

Case 1: $(\rho(L_{-\beta})(A_{\beta}))(\rho(L_{-\gamma})(A_{\gamma}))\subset A_{0,[\alpha]}.$
 In fact, for any $v_{-\beta}\in L_{-\beta}, a_{\beta}\in A_{\beta},v_{-\gamma}\in L_{-\gamma}, a_{\gamma}\in A_{\gamma}$,
we get $\rho(L_{-\beta})\in Der_{\phi}(A)$.
By Definition 2.3, we have
\begin{eqnarray*}
&&(\rho(v_{-\beta})(a_{\beta}))(\rho(v_{-\gamma})(a_{\gamma}))\\
&=&\rho(v_{-\beta})(a_{\beta}\phi^{-1}(\rho(v_{-\gamma})(a_{\gamma})))-
   \phi(a_{\beta})\rho(v_{-\beta})(\phi^{-1}(\rho(v_{-\gamma})(a_{\gamma}))).
\end{eqnarray*}
By Lemma 2.11-(5), $\rho(v_{-\gamma})(a_{\gamma})\in A_0$.
Since $A_\alpha$ is $\phi$-invariant, we have $\phi^{-1}(\rho(v_{-\gamma})(a_{\gamma}))\in A_0$
 and $a_{\beta}\phi^{-1}(\rho(v_{-\gamma})(a_{\gamma}))\in A_\beta$ by Lemma 2.11-(3).
 So $\rho(v_{-\beta})(a_{\beta}\phi^{-1}(\rho(v_{-\gamma})(a_{\gamma})))\in\rho(L_{-\beta})A_\beta$.
Similarly, one may check that  $\phi(a_{\beta})\rho(v_{-\beta})(\phi^{-1}(\rho(v_{-\gamma})(a_{\gamma})))\in A_\beta A_{-\beta}.$
Therefore, we get  $(\rho(L_{-\beta})(A_{\beta}))(\rho(L_{-\gamma})(A_{\gamma}))\subset A_{0,[\alpha]},$ as desired.

Case 2: $(\rho(L_{-\beta})(A_{\beta}))(A_{-\gamma}A_{\gamma})\subset A_{0,[\alpha]}.$
 In fact, for any $v_{-\beta}\in L_{-\beta}, a_{\beta}\in A_{\beta},v_{-\gamma}\in L_{-\gamma}, a_{\gamma}\in A_{\gamma}$,
since  $\rho(L_{-\beta})\in Der_{\phi}(A)$, we have
\begin{eqnarray*}
&&(\rho(v_{-\beta})(a_{\beta}))(a_{-\gamma}a_{\gamma})\\
&=&\rho(v_{-\beta})(a_{\beta}\phi^{-1}(a_{-\gamma}a_{\gamma}))-
   \phi(a_{\beta})\rho(v_{-\beta})(\phi^{-1}(a_{-\gamma}a_{\gamma})).
\end{eqnarray*}
By Lemma 2.11-(5), we get $\phi^{-1}(a_{-\gamma}a_{\gamma})\in A_0$
 and $a_{\beta}\phi^{-1}(a_{-\gamma}a_{\gamma})\in A_\beta.$
 Thus we have $\rho(v_{-\beta})(a_{\beta}\phi^{-1}(a_{-\gamma}a_{\gamma}))\in\rho(L_{-\beta})A_\beta$.
Similarly, one may check that  $\phi(a_{\beta})\rho(v_{-\beta})(\phi^{-1}(a_{-\gamma}a_{\gamma}))\\\in A_\beta A_{-\beta}.$
Therefore, we get  $(\rho(L_{-\beta})(A_{\beta}))(A_{-\gamma}A_{\gamma})\subset A_{0,[\alpha]},$  as desired.

Case 3: $(A_{-\beta}A_{\beta})(\rho(L_{-\gamma})(A_{\gamma}))\subset A_{0,[\alpha]}.$
In fact, by the commutativity of $A$, we have
$$(A_{-\beta}A_{\beta})(\rho(L_{-\gamma})(A_{\gamma}))=(\rho(L_{-\gamma})(A_{\gamma}))(A_{-\beta}A_{\beta})\subset A_{0,[\alpha]}.$$

Case 4: $(A_{-\beta}A_{\beta})(A_{-\gamma},A_{\gamma})\subset A_{0,[\alpha]}.$
In fact, if $\beta+\gamma=0$, that is, $\gamma=-\beta$, by the commutativity of $A$, we have
$$(A_{-\beta}A_{\beta})(A_{-\gamma}A_{\gamma})=A_{-\beta}(A_{\beta}A_{\beta}A_{-\beta})\subset A_{-\beta}A_{\beta}\subset A_{0,[\alpha]}.$$
If $\beta+\gamma\neq 0$, then
$$(A_{-\beta}A_{\beta})(A_{-\gamma},A_{\gamma})=(A_{-\beta}A_{-\gamma})(A_{\beta}A_{\gamma})\subset A_{-\beta-\gamma}A_{\beta+\gamma}\subset A_{0,[\alpha]}.$$

From four cases above, we have  $A_{0,[\alpha]}A_{0,[\alpha]}\subset A_{0,[\alpha]}\subset \mathcal{A}_{[\alpha]}$.
Similarly, one may check that
$A_{0,[\alpha]}A_{[\alpha]}\subset A_{[\alpha]}, A_{[\alpha]}A_{[\alpha]}\subset \mathcal{A}_{[\alpha]}.$
Therefore, we have $\mathcal{A}_{[\alpha]}\mathcal{A}_{[\alpha]}\subset \mathcal{A}_{[\alpha]}$.
The proof is finished.
$\hfill \Box$
 \medskip
%
%
%

\noindent{\bf Theorem 4.5.}
Let $A$ be a commutative and associative algebra associated to a Hom-Lie Rinehart algebra $L$.
Then the following assertions hold.

(1) For any $[\alpha]\in\Lambda/\approx$, the linear space
$\mathcal{A}_{[\alpha]}=A_{0,[\alpha]}\oplus A_{[\alpha]}$
of $A$ associated to $[\alpha]$ is an ideal of $A$.

(2) If $A$ is simple, then all weights of $\Lambda$ are connected. Furthermore,
$$A_{0}=(\sum_{-\alpha\in \Gamma,\alpha\in\Lambda}\rho(L_{-\alpha})(A_{\alpha})
                    +(\sum_{\alpha\in \Lambda}A_{-\alpha}A_{\alpha}).$$

\noindent{\bf Proof.}
(1) It is sufficient to show that $\mathcal{A}_{[\alpha]}A\subset \mathcal{A}_{[\alpha]}$ since $A$ is commutative.
In fact, by Proposition 4.4, we have
\begin{eqnarray*}
\mathcal{A}_{[\alpha]}A
=\mathcal{A}_{[\alpha]}(A_0\oplus(\bigoplus_{\beta\in[\alpha]}A_\beta)\oplus(\bigoplus_{\gamma\notin[\alpha]}A_\gamma))
\subset \mathcal{A}_{[\alpha]},
\end{eqnarray*}
as desired.

(2) Since $A$ is simple, then $\mathcal{A}_{[\alpha]}=A$ for any $\alpha\in\Lambda$.
So $[\alpha]=\Lambda$ and therefore
$A_{0}=(\sum_{-\alpha\in \Gamma,\alpha\in\Lambda}\rho(L_{-\alpha})(A_{\alpha})
                    +(\sum_{\alpha\in \Lambda}A_{-\alpha},A_{\alpha}).$
$\hfill \Box$
 \medskip

\noindent{\bf Theorem 4.6.}
Let $A$ be a commutative and associative algebra associated to a Hom-Lie Rinehart algebra $L$.
Then
\begin{eqnarray*}
A=V+\sum_{[\alpha]\in\Lambda/\approx}\mathcal{A}_{[\alpha]},
\end{eqnarray*}
where $V$ is a linear complement in $A_0$ of
$(\sum_{-\alpha\in \Gamma,\alpha\in\Lambda}\rho(L_{-\alpha})(A_{\alpha})
                    +(\sum_{\alpha\in \Lambda}A_{-\alpha},A_{\alpha})$
and any $A_{[\alpha]}$ is one of the ideals of $A$ described in Theorem 4.5-(1),
satisfying $\mathcal{A}_{[\alpha]}\mathcal{A}_{[\beta]}=0$, whenever $[\alpha]\neq[\beta]$.
\medskip

\noindent{\bf Proof.}
Straightforward from Theorem 4.5 and Proposition 4.4. $\hfill \Box$
 \medskip

Let us denote by $Z(A)$ the center of $A$, that is,
$Z(A):=\{a\in L| aA=0\}.$
\medskip

\noindent{\bf Corollary 4.7.}
Let $(L,A)$ be a  Hom-Lie Rinehart algebra.
If $Z(A)=0$ and
$$A_{0}=(\sum_{-\alpha\in \Gamma,\alpha\in\Lambda}\rho(L_{-\alpha})(A_{\alpha})
                    +(\sum_{\alpha\in \Lambda}A_{-\alpha}A_{\alpha}),$$
then $A$ is the direct sum of the ideals given in Theorem 4.5, that is,
    \begin{eqnarray*}
A=\sum_{[\alpha]\in\Lambda/\approx}\mathcal{A}_{[\alpha]},
\end{eqnarray*}
  satisfying $\mathcal{A}_{[\alpha]}\mathcal{A}_{[\beta]}=0$, whenever $[\alpha]\neq[\beta]$.
\medskip

\noindent{\bf Proof.}
This can be proved  analogously to Corollary 3.8 in \cite{Aragon2015}. $\hfill \Box$
 \medskip

In the following, we will discuss the relation between the decompositions of $L$ and $A$ of a Hom-Lie Rinehart algebra $(L,A)$.
\medskip

\noindent{\bf Definition 4.8.}
A split Lie-Rinehart algebra $(L,A)$ is \emph{tight} if $Z(L)=0,Z(A)=0,AA=A,AL=L$ and
\begin{eqnarray*}
   H=(\sum_{\gamma\in \Gamma,-\gamma\in\Lambda}A_{-\gamma}L_{\gamma})+(\sum_{\gamma\in \Gamma}[L_{-\gamma},L_{\gamma}]),
~A_{0}=(\sum_{-\alpha\in \Gamma,\alpha\in\Lambda}\rho(L_{-\alpha})(A_{\alpha}))
                    +(\sum_{\alpha\in \Lambda}A_{-\alpha}A_{\alpha}).
\end{eqnarray*}

\noindent{\bf Remark 4.9.}
Let $(L,A)$ be a tight split regular Hom-Lie-Rinehart algebra, then
    \begin{eqnarray*}
L=\sum_{[\gamma]\in\Gamma/\sim}I_{[\gamma]},~A=\sum_{[\alpha]\in\Lambda/\approx}\mathcal{A}_{[\alpha]},
\end{eqnarray*}
with any $I_{[\gamma]}$ an ideal of $L$ verifying $[I_{[\gamma]},I_{[\delta]}]=0$ if $ [\gamma]\neq[\delta] $
and any $\mathcal{A}_{[\alpha]}$ an ideal of $A$
satisfying  $\mathcal{A}_{[\alpha]}\mathcal{A}_{[\beta]}=0$ if $ [\alpha]\neq[\beta].$
\medskip

\noindent{\bf Proposition 4.10.}
Let $(L,A)$ be a tight split regular Hom-Lie Rinehart algebra, then for any $[\gamma]\in\Gamma/\sim$
there exists a unique $[\alpha]\in\Lambda/\approx$ such that $\mathcal{A}_{[\alpha]}I_{[\gamma]}=0$.
 \medskip

\noindent{\bf Proof.}
Similar to Proposition 4.2 in \cite{Aragon2015}.
$\hfill \Box$
 \medskip

\noindent{\bf Theorem 4.11.}
Let $(L,A)$ be a tight split regular Hom-Lie Rinehart algebra, then
    \begin{eqnarray*}
L=\sum_{i\in\Gamma/ I}L_{i},~A=\sum_{j\in J}A_{j},
\end{eqnarray*}
with any $L_i$ a nonzero ideal of $L$ and any $A_j$ a nonzero ideal of $A$.
Furthermore, for any $i\in I$ there exists a unique $\tilde{i}\in J$ such that $A_{\tilde{i}}L_i=0$.
 \medskip

\noindent{\bf Proof.}
Straightforward from Proposition 4.10.$\hfill \Box$

\section{The simple components}
\def\theequation{\arabic{section}. \arabic{equation}}
\setcounter{equation} {0}

In this section we focus on the simplicity of split regular Hom-Lie Rinehart algebra $(L,A)$ by centering our attention in those of maximal length.
From now on we always assume that $\Gamma$ and $\Lambda$ are symmetric in the sense that $\Gamma=-\Gamma$ and $\Lambda=-\Lambda$.
 \medskip

\noindent{\bf Lemma 5.1.}
 Let $(L,A)$ be a split regular Hom-Lie Rinehart algebra and $I$ an ideal of $L$.
Then $I=(I\cap H)\oplus(I\cap \bigoplus_{\gamma\in\Gamma}L_\gamma)$.
  \medskip

\noindent{\bf Proof.}
 Since $(L,A)$ is split, we get $L=H\oplus(\bigoplus_{\gamma\in\Gamma}L_\gamma)$.
 By the assumption that $I$ is an ideal of $L$, it is clear that $I$ is a submodule of $L$.
 Thus $I$ is a  weight module and therefore
$I=(I\cap H)\oplus(I\cap \bigoplus_{\gamma\in\Gamma}L_\gamma)$.$\hfill \Box$
  \medskip

 \noindent{\bf Lemma 5.2.}
 Let $(L,A)$ be a split regular Hom-Lie Rinehart algebra with $Z(L)=0$ and $I$ an ideal of $L$.
If $I \subset H$, then $I=\{0\}$.
   \medskip

\noindent{\bf Proof.}
Since  $I \subset H$, $[I,H]\subset [H,H]=0$.
It follows that $[I,L]=[I,\bigoplus_{\gamma\in\Gamma}L_\gamma]\subset H\cap (\bigoplus_{\gamma\in\Gamma}L_\gamma)=0$.
So  $I \subset Z(L)=0$.
$\hfill \Box$
  \medskip

\noindent{\bf Definition 5.3.}
A split regular Hom-Lie Rinehart algebra
$(L,A)$ is called \emph{root-multiplicative} if  the following conditions hold for any $\gamma,\delta\in\Gamma$ and $\alpha,\beta\in\Lambda$:

(1) If $\gamma\psi^{-1}+\delta\psi^{-1}\in\Gamma$, then $[L_{\gamma},L_{\delta}]\neq 0.$

(2) If $\alpha+\gamma\in\Gamma$, then $A_\alpha L_{\gamma}\neq 0.$

(3) If $\alpha+\beta\in\Lambda$, then $A_\alpha A_\beta\neq 0.$

 \medskip

\noindent{\bf Definition 5.4.}
A split regular Hom-Lie Rinehart algebra $(L,A)$ is called \emph{of maximal length} if
$dim L_{\gamma}=dim A_\alpha=1$  for any $\gamma\in\Gamma$ and $\alpha\in\Lambda$.
 \medskip


\noindent{\bf Theorem 5.5.}
Let $(L,A)$ be a tight split regular Hom-Lie Rinehart algebra of  maximal length and root-multiplicative.
If  all the nonzero roots of $L$ are connected, then either $L$ is simple or
$L=I\oplus I'$ with simple ideals $I,I'$.
 \medskip

\noindent{\bf Proof}
By Lemma 5.1, for any nonzero ideal $I$ of $L$, $I$ can be written as
  $$I=(I\cap H)\oplus(\bigoplus_{\gamma\in\Gamma_{I}}I_\gamma),$$
where $I_\gamma=I\cap L_\gamma$ and $\Gamma_{I}=\{\gamma\in\Gamma|I_\gamma\neq 0\}$ for at least one $\gamma\in\Gamma$.

In the first case, we assume that there exists $\gamma\in\Gamma$ satisfying $0\neq L_\gamma\subset I$.
Since $\psi(I)=I$ and $\psi(L_\gamma)\subset L_{\gamma\psi^{-1}},\psi^{-1}(L_\gamma)\subset L_{\gamma\psi}$,
we can assert that if $\gamma\in\Gamma_I$ then $\{\gamma\psi^z|z\in Z\}\subset \Gamma_I$.
That is,
$
\{L_{\gamma\psi^z}|z\in Z\}\subset  I.
$

Now, we take any $\delta\in\Gamma$ satisfying $\delta\notin \{\pm\gamma\psi^z\}$.
Since the root $\Gamma$ is connected to $\delta$, there is a connection
$\{\zeta_1,\zeta_2,\cdots,\zeta_n\}\subset\Lambda\cup\Gamma$,
     with $n\geq 2$, such that

  $~~~~\zeta_1=\gamma\psi^{k}$ for some $k\in Z$,

  $~~~~\zeta_1\psi^{-1}+\zeta_2\psi^{-1}\in\Gamma$,

    $~~~~\zeta_1\psi^{-2}+\zeta_2\psi^{-2}+\zeta_3\psi^{-1}\in\Gamma$,

    $~~~~\zeta_1\psi^{-3}+\zeta_2\psi^{-3}+\zeta_3\psi^{-2}+\zeta_4\psi^{-1}\in\Gamma$,

    $~~~~~~~~~\cdots\cdots\cdots$

    $~~~~\zeta_1\psi^{-i}+\zeta_2\psi^{-i}+\zeta_3\psi^{-i+1}+\cdots+\zeta_{i+1}\psi^{-1}\in\Gamma$,

    $~~~~~~~~~\cdots\cdots\cdots$

    $~~~~\zeta_1\psi^{-n+2}+\zeta_2\psi^{-n+2}+\zeta_3\psi^{-n+3}+\cdots+\zeta_{n-1}\psi^{-1}\in\Gamma$.

  $~~~~~\zeta_1\psi^{-n+1}+\zeta_2\psi^{-n+1}+\zeta_3\psi^{-n+2}+\cdots+\zeta_{n}\psi^{-1}=\varepsilon\delta\psi^{-m}$
   for some $m\in Z$ and $\varepsilon\in\{\pm 1\}.$

Taking into account that $\zeta_1\in \Gamma_{I}$, we have that if $\zeta_2\in \Gamma$,
the root-multiplicativity and the maximal length of $L$ allow us to assert
$0\neq [L_{\zeta_1},L_{\zeta_2}]=L_{\zeta_1\psi^{-1}+\zeta_2\psi^{-1}}.$
Since $0\neq L_{\zeta_1}\subset I$, we have
$$0\neq L_{\zeta_1\psi^{-1}+\zeta_2\psi^{-1}}\subset I.$$
A similar argument applied to $\zeta_1\psi^{-1}+\zeta_2\psi^{-1}\in\Gamma$, $\zeta_3\in\Lambda\cup\Gamma$ and
$$
(\zeta_1\psi^{-1}+\zeta_2\psi^{-1})\psi^{-1}+\zeta_3\psi^{-1}
=\zeta_1\psi^{-2}+\zeta_2\psi^{-2}+\zeta_3\psi^{-1}\in\Gamma,
$$
we have
$$0\neq L_{\zeta_1\psi^{-2}+\zeta_2\psi^{-2}+\zeta_3\psi^{-1}}\subset I.$$
 We can follow this process with the connection $\{\zeta_1,\zeta_2,\cdots,\zeta_n\}$ to get
$$0\neq L_{\zeta_1\psi^{-n+1}+\zeta_2\psi^{-n+1}+\zeta_3\psi^{-n+2}+\cdots+\zeta_{n}\psi^{-1}}\subset I.$$
That is, $L_{\delta\psi^{-m}}\subset I$ or $L_{-\delta\psi^{-m}}\subset I$.

Since $-\gamma\in\Gamma$, it is easy to see that$\{-\zeta_1,-\zeta_2,\cdots,-\zeta_n\}$
 is a connection from $-\gamma$ to $\delta$ satisfying
$$-\zeta_1\psi^{-n+1}-\zeta_2\psi^{-n+1}-\zeta_3\psi^{-n+2}-\cdots-\zeta_{n}\psi^{-1}=-\varepsilon\delta\psi^{-m}.$$
By similar discussions above, one may check that $0\neq L_{-\varepsilon\delta\psi^{-m}}\subset I$.
Therefore, $\Gamma_I=\Gamma$ and $H\subset I$.
If we consider now any $\gamma\in\Gamma$, since $L_\gamma=[H,L_\gamma]$  and $H\subset I$, by the maximal length of $L$,
 we have $L_\gamma  \subset I$ and so $I=L$. That is, $L$ is simple.

In the second case, suppose that for any $\gamma\in\Gamma_I$ we have that $-\gamma\notin\Gamma_I$.
Then we have
  \begin{eqnarray}
\Gamma=\Gamma_I\dot{\cup} -\Gamma_I,
\end{eqnarray}
where $-\Gamma_I:=\{-\gamma|\gamma\in\Gamma_I\}$.
Define
  \begin{eqnarray}
I':=(\sum_{-\gamma\in -\Gamma_{I},\gamma\in\Lambda}A_{\gamma}L_{-\gamma})\oplus(\bigoplus_{-\gamma\in -\Gamma_{I}}L_{-\gamma}).
\end{eqnarray}

First, we claim  that  $I'$ is a Hom-Lie ideal of $L$.
In fact, by Eq. (3.5), $\psi(L_{-\gamma})\subset L_{-\gamma\psi^{-1}}$, $-\gamma\psi^{-1}\in -\Gamma_{I}$.
By Lemma 2.11-(4), $\psi(A_{\gamma}L_{-\gamma})\subset \psi(L_0)\subset L_0$  if $A_{\gamma}L_{-\gamma}\neq 0$ (otherwise is trivial).
So  $\psi(I')\subset I'$.

Since $A_{\gamma}L_{-\gamma}\subset L_0$, by Eq. (5.2), we have
  \begin{eqnarray}
[L,I']&=&[H\oplus(\bigoplus_{\delta\in\Gamma}L_\delta),
(\sum_{-\gamma\in -\Gamma_{I},\gamma\in\Lambda}A_{\gamma}L_{-\gamma})\oplus(\bigoplus_{-\gamma\in -\Gamma_{I}}L_{-\gamma})]\subset\nonumber\\
&&[\bigoplus_{\delta\in\Gamma}L_\delta,(\sum_{-\gamma\in -\Gamma_{I},\gamma\in\Lambda}A_{\gamma}L_{-\gamma})]
        +[\bigoplus_{\delta\in\Gamma}L_\delta,\bigoplus_{-\gamma\in -\Gamma_{I}}L_{-\gamma}]
        +\sum_{-\gamma\in -\Gamma_{I}}L_{-\gamma}.~~
\end{eqnarray}

For the expression $[\bigoplus_{\delta\in\Gamma}L_\delta,(\sum_{-\gamma\in -\Gamma_{I},\gamma\in\Lambda}A_{\gamma}L_{-\gamma})]$ in Eq. (5.3).
If some $[L_\delta,A_{\gamma}L_{-\gamma}]\neq 0,$
we have that in case $\delta=-\gamma$,
$[L_{-\gamma},A_{\gamma}L_{-\gamma}]\subset  L_{-\gamma\psi^{-1}}\subset I'$,
and in case $\delta=\gamma$, since $I$ is a Hom-Lie ideal of $L$, $-\gamma\notin \Gamma_{I}$ implies
$[L_{-\gamma},A_{-\gamma}L_{\gamma}]=0$, by the  maximal length of $L$ and the symmetry of $\Gamma$,
we have $[L_{\gamma},A_{\gamma}L_{-\gamma}]=0$.
Suppose $\delta\notin \{\gamma,-\gamma\}.$
 By Definition 2.3,
 \begin{eqnarray*}
[L_\delta,A_{\gamma}L_{-\gamma}]
\subset \phi(A_{\gamma})[L_\delta, L_{-\gamma}]+\rho(L_\delta)(A_{\gamma})L_{-\gamma}.
\end{eqnarray*}
Since $(L,A)$ is regular, $\phi(A_{\gamma})\subset A_{\gamma}.$
As $[L_\delta,A_{\gamma}L_{-\gamma}]\neq 0,$  we get
$A_{\gamma}[L_\delta, L_{-\gamma}]\neq 0$ or $\rho(L_\delta)(A_{\gamma})L_{-\gamma}\neq 0$.
By the  maximal length of $L$, either $A_{\gamma}[L_\delta, L_{-\gamma}]=L_{\gamma+(\delta-\gamma)\psi^{-1}}$
or $\rho(L_\delta)(A_{\gamma})L_{-\gamma}=L_{\delta}$.
In both cases, since $\gamma\in\Gamma_I$, by the root-multiplicativity of $L$,
we have $L_{-\delta}\subset I$ and therefore $-\delta\in\Gamma_I$. That is, $L_{\delta}\subset I'$.
So $[\bigoplus_{\delta\in\Gamma}L_\delta,(\sum_{-\gamma\in -\Gamma_{I},\gamma\in\Lambda}A_{\gamma}L_{-\gamma})]\subset I'$.

For the expression $[\bigoplus_{\delta\in\Gamma}L_\delta,\bigoplus_{-\gamma\in -\Gamma_{I}}L_{-\gamma}]$ in Eq. (5.3),
If some $[L_\delta,L_{-\gamma}]\neq 0,$ then $[L_\delta,L_{-\gamma}]=L_{(\delta-\gamma)\psi^{-1}}$.
On the one hand, let $\delta-\gamma\neq 0$.
Since  $\gamma\in\Gamma_I$, by the root-multiplicativity of $L$,
we have $[L_\gamma,L_{-\delta}]=L_{(\gamma-\delta)\psi^{-1}}\subset I$.
So $(\delta-\gamma)\psi^{-1}\in\Gamma_I$ and therefore $L_{(\delta-\gamma)\psi^{-1}}\subset I'$.
On the one hand, let $\delta-\gamma=0$. Suppose $[L_\gamma,L_{-\gamma}]\neq 0$, since $\gamma\in\Gamma_I$,
we get $[L_\gamma,L_{-\gamma}]\subset I$.
Thus $L_{-\gamma}=[[L_\gamma,L_{-\gamma}],L_{-\gamma\psi}]\subset I$.
According to the discussion above, $\gamma,-\gamma\in\Gamma_I$,
 a contradiction with Eq. (5.1).
So $[\bigoplus_{\delta\in\Gamma}L_\delta,\bigoplus_{-\gamma\in -\Gamma_{I}}L_{-\gamma}]\subset I'$.

Second, we claim that $\rho(I')(A)L\subset I'$.
In fact, by Definition 2.3, we have
 \begin{eqnarray*}
\rho(I')(A)L\subset [I',AL]+A[I',L]
 \end{eqnarray*}
Since $I'$ is a Hom-Lie ideal of $L$, we get $[I',AL]\subset I', [I',L]\subset I'$.
So it is sufficient to verify that $AI'\subset I'$. For this, we calculate
 \begin{eqnarray}
AI'&=&(A_0\oplus(\bigoplus_{\alpha\in\Lambda}A_\alpha))
(\sum_{-\gamma\in -\Gamma_{I},\gamma\in\Lambda}A_{\gamma}L_{-\gamma})\oplus(\bigoplus_{-\gamma\in -\Gamma_{I}}L_{-\gamma})\subset\nonumber\\
&&I'+(\bigoplus_{\alpha\in\Lambda}A_\alpha)(\sum_{-\gamma\in -\Gamma_{I},\gamma\in\Lambda}A_{\gamma}L_{-\gamma})
+(\bigoplus_{\alpha\in\Lambda}A_\alpha)(\bigoplus_{-\gamma\in -\Gamma_{I}}L_{-\gamma}).
 \end{eqnarray}

 For the expression $(\bigoplus_{\alpha\in\Lambda}A_\alpha)(\sum_{-\gamma\in -\Gamma_{I},\gamma\in\Lambda}A_{\gamma}L_{-\gamma})$ in Eq. (5.4),
 if some $A_\alpha(A_{\gamma}L_{-\gamma})\neq 0$, we have that in case $\alpha=-\gamma$, clearly
 $A_\alpha(A_{\gamma}L_{-\gamma})=A_{-\gamma}(A_{\gamma}L_{-\gamma})\subset L_{-\gamma}\subset I'$.
 In case of  $\alpha=\gamma$, since $-\gamma\notin \Gamma_I$, we get $A_{-\gamma}(A_{-\gamma}L_{\gamma})=0$.
 By the   the  maximal length of $L$, we have $A_\alpha(A_{\gamma}L_{-\gamma})=A_{\gamma}(A_{\gamma}L_{-\gamma})=0$.
 Suppose that $\alpha\notin\{\gamma,-\gamma\}$,
by the   the  maximal length of $L$, we have $A_\alpha(A_{\gamma}L_{-\gamma})=L_\alpha$.
Since $\gamma\in\Gamma_I$, by the root-multiplicativity of $L$, we have $L_{-\gamma})\subset I$, that is, $-\alpha\in\Gamma_I$.
 So $\alpha\in-\Gamma_I$ and $L_\alpha\subset I'$.
 Thus $(\bigoplus_{\alpha\in\Lambda}A_\alpha)(\sum_{-\gamma\in -\Gamma_{I},\gamma\in\Lambda}A_{\gamma}L_{-\gamma})\subset I'$.

  For the expression $(\bigoplus_{\alpha\in\Lambda}A_\alpha)(\bigoplus_{-\gamma\in -\Gamma_{I}}L_{-\gamma})$ in Eq. (5.4),
if some $A_\alpha L_{-\gamma}\neq 0$, in case $\alpha-\gamma\in \Gamma_{I}$,
by the root-multiplicativity of $L$, we have $A_{-\alpha} L_{\gamma}\neq 0$.
Again by the  maximal length of $L$, we have $A_{-\alpha} L_{\gamma}=L_{-\alpha+\gamma}$.
So $ -\alpha+\gamma\in \Gamma_{I}$, a contradiction. Thus $ \alpha-\gamma\in -\Gamma_{I}$ and therefore
$(\bigoplus_{\alpha\in\Lambda}A_\alpha)(\bigoplus_{-\gamma\in -\Gamma_{I}}L_{-\gamma})\subset I'$.

By the discussion above, we have shown that $\rho(I')(A)L\subset I'$ and therefore $I'$ is an ideal of $(L,A)$.

Finally, we will verify that $L=I\oplus I'$ with simple $I,I'$.
By similar discussions in the proof  of the simplicity of $L$ above, it is not hard to check that $I,I'$ are simple Hom-Lie ideals.
Since $[I',I]=0$, by the commutativity of  $H$, we get $\sum_{\gamma\in\Gamma}[L_{\gamma},L_{-\gamma}]=0$,
so $H$ must has the form
  \begin{eqnarray}
H=(\sum_{\gamma\in \Gamma_{I},-\gamma\in\Lambda}A_{\gamma}L_{-\gamma})\oplus(\sum_{-\gamma\in -\Gamma_{I},\gamma\in\Lambda}A_{\gamma}L_{-\gamma}).
\end{eqnarray}
In order to show that the sum in Eq. (5.5) is direct, we take any
$h\in(\sum_{\gamma\in \Gamma_{I},-\gamma\in\Lambda}A_{\gamma}L_{-\gamma})\oplus(\sum_{-\gamma\in -\Gamma_{I},\gamma\in\Lambda}A_{\gamma}L_{-\gamma}).$
Suppose $h\neq 0$, then $h\notin Z(L)$.
Since $L$ is split, there is $v_\delta\in L_\delta, \delta\in \Gamma$ satisfying
$[h,v_\delta]=\delta(h)\psi(v_\delta)\neq 0$.
By Eq. (3.5), $0\neq\delta(h)\psi(v_\delta)\in L_{\delta\psi^{-1}}$.
While $L_{\delta\psi^{-1}}\subset I\cap I'=0$, a contradiction.
So $h=0$, as required.
The proof is completed.
$\hfill \Box$
  \medskip

\noindent{\bf Theorem 5.6.}
Let $(L,A)$ be a tight split regular Hom-Lie Rinehart algebra of  maximal length and root-multiplicative.
If  all the nonzero weights of $A$ are connected, then either $L$ is simple or
$A=J\oplus J'$ with simple ideals $J,J'$.
 \medskip

\noindent{\bf Proof.}
 This can be proved completely analogously to Theorem 5.5.
$\hfill \Box$
  \medskip

\noindent{\bf Theorem 5.7.}
Let $(L,A)$ be a tight split regular Hom-Lie Rinehart algebra of  maximal length and root-multiplicative.
Assume that  the roots and weights systems are both  symmetric such that $\Gamma$ have all its
nonzero roots connected and $\Lambda$ have all its nonzero weights connected.
Then
\begin{eqnarray*}
L=\bigoplus_{i\in I}L_i,~~A=\bigoplus_{j\in J}A_j,
\end{eqnarray*}
where any $L_i$ is a simple ideal of $L$ having all of its nonzero roots connected satisfying $[L_i,L_{i'}]=0$ for any $i'\in I$ with $i\neq i'$,
and any $A_j$ is a simple ideal of $A$ satisfying $A_j A_{j'}=0$ for any $j'\in J$ such that $j'\neq j$.

Furthermore, for any $i\in I$ there exists a unique $\overline{i}\in J$ such that $A_{\overline{i}}L_i\neq 0$.
We also have that any $L_i$ is a split regular Hom-Lie Rinehart algebra over $A_{\overline{i}}$.
  \medskip

\noindent{\bf Proof.}
  It is analogous to Theorem 5.5 in \cite{Albuquerque}.
$\hfill \Box$

\begin{center}
 {\bf ACKNOWLEDGEMENT}
 \end{center}

  The paper is  supported by
   the Anhui Provincial Natural Science Foundation (Nos. 1908085MA03 and 1808085MA14),
   the NSF of China (Nos. 11761017 and 11801304), and the Youth Project for Natural Science Foundation of Guizhou provincial department of education (No. KY[2018]155).

\renewcommand{\refname}{REFERENCES}

\end{document}